\documentclass{amsart}
\usepackage{amssymb}
\newtheorem{theorem}{Theorem}
\newtheorem{lemma}[theorem]{Lemma}

\title{Beyond Undecidable} 
\author{Paola Cattabriga} 
\address{Universit\`a di Bologna,   ITALY} 
\email{co14099@iperbole.bologna.it} 
\begin{document}
\bigskip
\begin{abstract}  The predicate complementary  to the well-known G\"{o}del's 
provability predicate is defined.
From its recursiveness new consequences concerning the 
incompleteness argumentation are drawn and extended to new results of 
consistency,
completeness and decidability with regard to Peano Arithmetic and the
first order predicate calculus.  \end{abstract}
\thanks{Appeared in Proceedings of The International 
Conference on Artificial Intelligence, IC-AI'2000, H. R. Arabnia (ed.), 
Las Vegas, Nevada, USA, June 26-29, 2000, Volume III,  CSREA Press, pp. 1475-1481}
 \maketitle
 
\bigskip\bigskip

\section*{Introduction}
Of all the remarkable logical achievements of the twentieth 
century perhaps the most outstanding is the celebrated  G\"{o}del
incompleteness argumentation of 1931 \cite{godel1,godel2}.
In contrast to Hilbert's program called for embodying classical
mathematics in a formal system and proving that system
consistent by finitary methods \cite{hilbert2},
G\"{o}del paper showed that not even the first step could be carried 
out fully,  any formal system suitable for the arithmetic of 
integers was incomplete.

The present article, in the most 
absolute respect for the extraordinary contribution given by  
G\"{o}del to the logical inquiry, 
brings  G\"{o}del's achievement into question  by the definition 
of the refutability predicate.
As it is well-known self-reference plays a crucial role in 
G\"{o}del's incompleteness argumentation and the  methods of 
achieving self-referential statements  
is the so-called ``diagonalization". 
The refutability predicate, defined by arithmetization  
 as a number-theoretic statement, 
gives rise to new consequences properly regarding  
G\"{o}del's incompleteness argumentation and the method of 
diagonalization.   This article proposes  a revision based on the logical
investigation of the interactive links between provability and refutability
 predicates.
Originally devised by G\"{o}del in order to arithmetize 
metamathematical notions, G\"{o}del numbering turns out to be the key
of the problem in defining  refutability with the same recursive 
status as  provability. 
The inquiry comes up with a final solution for finitary methods and 
the related decision problem \cite{hilbert1}.

The paper is organized as follows. Firstly, in the following of 
this section, we introduce diagonalization and the famous 
incompleteness argumentation of G\"{o}del. 
Section \ref{R} presents two new primitive recursive predicates for 
refutability and the enucleation of some of their consequences, 
which represent the first main result of this paper: 
G\"{o}del's incompleteness argumentation is not a
theorem in Peano Arithmetic. Section \ref{C} shows that any 
formula of Peano Arithmetic is  proved if and only if it is not 
refuted, and extends this result to the accomplishment of consistency and 
completeness for Peano Arithmetic and then to the achievement of decidability for 
first order predicate calculus.  

\bigskip

\subsection*{Basic Setup}
We shall assume a first order theory which adequately formalizes 
Peano Arithmetic (see 
for example the system $S$, with all the necessary assumptions, in 
\cite{mendel} 116-175).
 Let us call it {\it PA}. 
As  is well known by means of the G\"{o}del numbering, each 
expression in {\it PA} can 
refer to itself.  Numerals, as usual,
are  defined recursively, $\overline{0}$ is $0$ and for any natural 
number $n$, $\overline{n+1}$
is $(\overline{n})'$ (where $'$ is the Successor function). For any 
expression $X$  we use $\ulcorner X \urcorner$ to denote the 
corresponding
G\"{o}del number of $X$. Let us define the G\"{o}del numbering as 
follows:
\begin{enumerate}
\item First assign different odd numbers to the primitive symbols of 
the language
of {\it PA}.
\item Let $X$ be a formal expression $X_{0},X_{1},\dots,X_{n}$, where 
each $X_i$, 
$0\leqslant i\leqslant n$, is
a primitive symbol of the language of  {\it PA}. Then
$$ \ulcorner X\urcorner\,=\,p_{0}^{\ulcorner 
X_{0}\urcorner}\centerdot 
p_{1}^{\ulcorner X_{1}\urcorner}\centerdot\ldots\centerdot 
p_{n}^{\ulcorner X_{n}\urcorner}$$
where $p_n$ is the $n$-th prime number and $p_{0}=2$.
\item Let $X$ be composed by the formal expressions 
$X_{0},X_{1},\dots,X_{n}$, then
$$ \ulcorner X\urcorner\,=\,p_{0}^{\ulcorner 
X_{0}\urcorner}\centerdot 
\ldots\centerdot p_{n}^{\ulcorner X_{n}\urcorner}.$$
\end{enumerate}
For any given  formula  $\phi(v)$ of {\it PA} we then have  its 
G\"{o}del number 
$n\,=\,\ulcorner\phi (v)\urcorner$. This number $n$ has a name in the 
language
 of {\it PA},
namely $\overline n$, and this name can be substituted back into 
$\phi(v)$. This
self-reference procedure is admitted by the so-called 
\emph{diagonalization lemma} as  follows.

\bigskip

\subsection*{Diagonalization}
For any formula $\phi$ with only the variable $v$ free there is a
sentence $\delta$ such that
$$
\vdash _{\mathit{PA}}\delta\, \iff \,\phi (\overline{\ulcorner \delta 
\urcorner}).
$$
The argumentation usually considered to be a proof is  the following.
We  define the function of substitution 
$sb(\ulcorner \phi (v)\urcorner,\overline{n})\,=\,\ulcorner \phi 
(\overline{n}) \urcorner $,
which gives us the G\"{o}del number of the result of replacing $v$ by 
the $n$-th numeral
in $\phi (v)$
 (see the corresponding $Sb(x^{v}_{y})$ and $Sb[x^{a}_{\chi (y)}]$ in 
\cite{godel1,godel2}).

Let $\phi (v)$ be given and let us call $\beta (v)$ the formula $\phi 
(sb(v,v))$. Let 
$m\,=\,\ulcorner \beta (v) \urcorner$ and $\delta\, =\, 
\beta(\overline{m})$. We shall
show that $\delta$ is the sentence we were looking for. To this 
purpose we notice that
in {\it PA} they hold the following equivalences
\begin{xalignat}{2} \notag
\vdash \delta &\iff  \beta(\overline{m}) && \text{by definition}\\ \notag
& \iff \phi (sb(\overline{m},\overline{m})) && \text{by definition}\\ \notag
& \iff \phi (sb(\ulcorner\overline{\beta (v)}\urcorner,\overline{m})) 
 && \text{since } m=\ulcorner \beta (v) \urcorner\\ \notag
& \iff \phi (\ulcorner\overline{\beta (\overline{m})}\urcorner)  
&&\text{definition of }sb\\ \notag
& \iff \phi (\overline{\ulcorner \delta\urcorner}) && \text{by definition.}
\end{xalignat}

\bigskip

\subsection*{G\"{o}del's Incompleteness}
We  present the version of  the so-called  \emph{G\"{o}del's  first 
incompleteness Theorem} 
as it is given in (\cite{mendel} 161-162), to which the reader can 
refers for
the definition of the concepts which are involved.

Let $\phi(v)$ be the formula $\forall x\,\neg Pf(x,v)$, hence
by diagonalization lemma we attain $$\vdash _{\mathit PA} \delta \iff \forall x\,\neg 
Pf(x,\overline{\ulcorner \delta\urcorner}).$$
G\"{o}del's   incompleteness argumentation asserts:
\begin{itemize}
\item[(a)] if {\it PA} is consistent, not $\vdash _{PA} \delta$,
\item[(b)] if {\it PA} is $\omega$-consistent, not $\vdash _{PA} 
\neg\delta$,
\end{itemize}
hence, if {\it PA} is $\omega$-consistent, $ \delta$ is an 
undecidable sentence of {\it PA}.

The proof is as follows. Let $q$ be the G\"{o}del number of $\delta$.
\begin{itemize}
\item[(a)] Assume $\vdash _{PA} \delta$. Let $r$ be the G\"{o}del 
number of a proof in {\it PA} of $\delta$. Then $\mathtt{Pf}$$(r,q)$. Hence, $\vdash 
_{PA} Pf(\overline{r},\overline{q})$, that is $\vdash _{PA} 
Pf(\overline{r},\overline{\ulcorner \delta\urcorner})$. We already 
have $\vdash _{PA} \delta \iff \forall x\,\neg Pf(x,\overline{\ulcorner \delta\urcorner})$. By 
Biconditional Elimination,
$\vdash _{PA} \forall x\,\neg Pf(x,\overline{\ulcorner \delta\urcorner})$. By Rule A4 
(Particularization Rule),
$\vdash _{PA} \neg Pf(\overline{r},\overline{\ulcorner 
\delta\urcorner})$. Therefore,
{\it PA} is inconsistent.
\item[(b)] Assume {\it PA} is $\omega$-consistent and  $\vdash _{PA} \neg\delta$. Since
$\vdash _{PA} \delta \iff \forall x\,\neg Pf(x,\overline{\ulcorner \delta\urcorner})$, 
Biconditional Elimination  yields 
$\vdash _{PA} \neg \forall x\,\neg Pf(x,\overline{\ulcorner \delta\urcorner})$ which abbreviates
to $(\ast)\vdash _{PA} \exists x\, Pf(x,\overline{\ulcorner \delta\urcorner})$. On the other
hand, since {\it PA} is $\omega$-consistent, {\it PA} is consistent. 
But, $\vdash _{PA} \neg \delta$. Hence, not $\vdash _{PA}  \delta$; that 
is, there
is no proof in {\it PA} of $\delta$. So $\mathtt{Pf}$$(n,q)$ is false 
for every natural number
$n$ and, therefore, $\vdash _{PA} \neg 
Pf(\overline{n},\overline{\ulcorner \delta\urcorner})$
for every natural number $n$.(Remember that $\overline{\ulcorner \delta\urcorner}$ is $\overline{q}$.)
By $\omega$-consistency, not $\vdash _{PA} \exists x\, Pf(x,\overline{\ulcorner \delta\urcorner})$,
contradicting $(\ast)$.
\end{itemize}

\bigskip

\section{Refutability}\label{R}
We are now ready to present the results with which this paper is 
concerned.
We shall construct two new predicates by  G\"{o}del numbering. 
The reader can refer to the arithmetization as defined by Mendelson;  
the new predicates must be considered as two last relations added 
to the functions and  relations  (1-26) presented 
in (\cite{mendel} 149-156)\footnote{We shall not reproduce entirely 
this 
long list of definitions which is already  well-known (see also 
\cite{godel1} 162-176).}.  
Let us start  recalling  some of the definitions  involved, precisely 
only
those we need.
    
$\mathtt{MP}$$(x,y,z)$: The expression  with  G\"{o}del number $z$ is 
a direct consequence
of the expressions with G\"{o}del numbers $x$ and $y$ by modus ponens,
$$y = 2^{3} * x * 2^{11} * z * 2^{5} \wedge\mathtt{Gd}
\mathrm (x)\wedge\mathtt{Gd}\mathrm (z). $$
  
$\mathtt{Gen}$$(x,y)$: The expression  with  G\"{o}del number $y$ 
comes from
 the expression with G\"{o}del number $x$  by the Generalization Rule,
 
\smallskip

\noindent $(\exists v)_{v<y} (\mathtt{EVbl}\mathrm (v) \wedge 
y = 2^{3} *  2^{3} *  2^{13} * v * 2^{5} * x * 2^{5} \wedge 
\mathtt{Gd}\mathrm (x)). $
  
\smallskip
  
$\mathtt{Ax}$$(y)$: $y$ is the   G\"{o}del number of an axiom of {\it 
PA}:
$$\mathtt{LAx}\mathrm (y)\vee\mathtt{PrAx}\mathrm (y). $$
  
$\mathtt{Neg}$$(v)$: the G\"{o}del number of $(\neg \alpha)$ if $v$ 
is the  G\"{o}del 
number of $ \alpha$:
$$ \mathtt{Neg}\mathrm (v) = 2^{3} * 2^{9} * v * 2^{5}.$$
  
$\mathtt{Prf}$$(x)$: $x$ is the G\"{o}del number of a proof in {\it 
PA}: 

\smallskip

\noindent $
\exists u_{u<x}\; \exists v_{v<x} \;\exists z_{z<x}\; \exists 
w_{w<x}\; ([x  =2^{w} \wedge \mathtt{Ax}\mathrm (w)] \vee $

\noindent $
[\mathtt{Prf}\mathrm (u) \wedge  \mathtt{Fml}\mathrm ((u)_{w})  \wedge  x  = u * 2^{v}  
\wedge \mathtt{Gen}\mathrm ((u)_{w},v)]\vee $ 

\noindent$
[\mathtt{Prf}\mathrm (u) \wedge \mathtt{Fml}\mathrm ((u)_{z}) \wedge  
\mathtt{Fml}\mathrm ((u)_{w})  \wedge   x  = u * 2^{v} \wedge 
 \mathtt{MP}\mathrm ((u)_{z},(u)_{w},v)]\vee $

\noindent $[\mathtt{Prf}\mathrm (u) \wedge  x  = u * 2^{v} \wedge 
\mathtt{Ax}\mathrm (v)].
$
  
\smallskip
  
$\mathtt{Pf}$$(x,v)$: $x$ is the G\"{o}del number of a proof in {\it 
PA} of the
 formula with G\"{o}del number $v$: $$ \mathtt{Prf}\mathrm (x) \wedge 
v = (x)_{\mathit{lh}
 \mathrm (x) \overset{\centerdot}{\text{--}} 1}.$$

By means of such definitions, we shall define two new 
predicates,  $\mathtt{Rf}$ and $\mathtt{Ref}$.

\smallskip
 
$\mathtt{Rf}\mathrm (x,v)$: $x$ is the G\"{o}del number of a proof in 
{\it PA} of the negation of the formula with G\"{o}del 
number $v$: $$ \mathtt{Pf} \mathrm (x,z) \wedge z = \mathtt{Neg}\mathrm (v). $$

\noindent In other terms  $\mathtt{Rf} \mathrm (x,v)$ states $x$ 
is the G\"{o}del number of a refutation in 
{\it PA} of the formula with G\"{o}del number $v $ \footnote{One can 
easily  see that $ \mathtt{Rf} \mathrm (x,v) $ is the same as

\noindent $\exists u_{u<x}\; \exists v_{v<x} \;\exists z_{z<x}\; \exists 
w_{w<x}\;\exists y_{y<x}$ 

\noindent $([x=2^{y} \wedge \mathtt{Ax}\mathrm (y) \wedge 
y=\mathtt{Neg}\mathrm (v)] \vee $

\noindent $[\mathtt{Prf}\mathrm (u) \wedge \mathtt{Fml}\mathrm ((u)_{w}) \wedge 
x = u * 2^{y} \wedge \mathtt{Gen}\mathrm ((u)_{w},y)\wedge 
y=\mathtt{Neg}\mathrm (v)]\vee $

\noindent $[\mathtt{Prf}\mathrm (u) \wedge \mathtt{Fml}\mathrm ((u)_{z}) \wedge 
\mathtt{Fml}\mathrm ((u)_{w})  \wedge x = u * 2^{y} \wedge 
\mathtt{MP}\mathrm ((u)_{z},(u)_{w},y) \wedge 
y=\mathtt{Neg}\mathrm (v)]\vee $

\noindent $[\mathtt{Prf}\mathrm (u) \wedge  x = u * 2^{y} \wedge 
\mathtt{Ax}\mathrm (y) \wedge y=\mathtt{Neg}\mathrm (v)].$}.

 $\mathtt{Rf}$ is \emph{primitive recursive}, as the 
relations 
obtained from primitive recursive relations by means of 
propositional connectives are also primitive recursive (\cite{mendel} 
137).
For its recursiveness  $\mathtt{Rf} \mathrm (x,v)$  is
\emph{expressible} in {\it PA} by a formula $Rf(x,v)$.

 \smallskip

$\mathtt{Ref}$$(x)$: $x$ is the G\"{o}del number of a refutation in 
{\it PA}: 
$$\mathtt{Prf}\mathrm (v) \wedge v = \mathtt{Neg}\mathrm (x)  $$ 
 
 \smallskip
In other terms  $\mathtt{Ref} \mathrm (x)$ states $x$ 
is the G\"{o}del number of a proof in {\it PA} of its negation. 
$\mathtt{Ref}$$(x)$ is \emph{primitive recursive}, as the relations 
obtained from primitive recursive relations by means of 
propositional connectives are also primitive recursive.
For its recursiveness $\mathtt{Ref}$$(x)$ is \emph{expressible} in 
{\it PA} by 
a formula $Ref(x)$.
\begin{lemma}\label{notboth}
For any natural number $n$ and for any formula $\alpha$  not both
$\mathtt{Rf}\mathrm (n,\ulcorner \alpha \urcorner)$ 
and $\mathtt{Pf}\mathrm (n,\ulcorner \alpha \urcorner)$.
\end{lemma}
\begin{proof}
Let us suppose to have both $\mathtt{Rf}\mathrm (n,\ulcorner \alpha 
\urcorner)$ and
$\mathtt{Pf}\mathrm (n,\ulcorner \alpha \urcorner)$.
We should have then
$\mathtt{Prf}\mathrm (n) \wedge \ulcorner \alpha \urcorner = 
(n)_{\mathit{lh}
 \mathrm (n) \overset{\centerdot}{\text{--}} 1}$ and $\mathtt{Pf} 
\mathrm (n,z) \wedge 
z = \mathtt{Neg}\mathrm (\ulcorner \alpha \urcorner),$ 
i.e. 
$\mathtt{Prf}\mathrm (n) \wedge \ulcorner \alpha \urcorner = 
(n)_{\mathit{lh}
 \mathrm (n) \overset{\centerdot}{\text{--}} 1}$
 and  $\mathtt{Prf} \mathrm (n) \wedge 
\mathtt{Neg}\mathrm (\ulcorner \alpha \urcorner) = (n)_{\mathit{lh}
 \mathrm (n) \overset{\centerdot}{\text{--}} 1}.$

By the definition of $\mathtt{Prf} \mathrm (x) $ this would mean to 
have

\noindent $
\exists u_{u<n}\; \exists v_{v<n} \;\exists z_{z<n}\; \exists 
w_{w<n}\; ([n=2^{w} \wedge \mathtt{Ax}\mathrm (w)] \vee $

\noindent $[\mathtt{Prf}\mathrm (u) \wedge \mathtt{Fml}\mathrm ((u)_{w}) \wedge 
n = u * 2^{v} \wedge \mathtt{Gen}\mathrm ((u)_{w},v)]\vee $

\noindent $[\mathtt{Prf}\mathrm (u) \wedge \mathtt{Fml}\mathrm ((u)_{z}) \wedge 
\mathtt{Fml}\mathrm ((u)_{w})  \wedge n = u * 2^{v} \wedge 
\mathtt{MP}\mathrm ((u)_{z},(u)_{w},v)]\vee $

\noindent $[\mathtt{Prf}\mathrm (u) \wedge  n = u * 2^{v} \wedge 
\mathtt{Ax}\mathrm (v)]
$ 
and both
$$\mathrm \ulcorner \alpha \urcorner = (n)_{\mathit{lh}
 \mathrm (n) \overset{\centerdot}{\text{--}} 1}
\:\:\text{ and }\:\:
\mathtt{Neg}\mathrm (\ulcorner \alpha \urcorner) = (n)_{\mathit{lh}
 \mathrm (n) \overset{\centerdot}{\text{--}} 1}
$$
and hence the four cases
\begin{enumerate}
	\item  $[n=2^{\ulcorner \alpha \urcorner} \wedge \mathtt{Ax}\mathrm 
(\ulcorner \alpha \urcorner)] 
\text{ and }$

\noindent $ [n=2^{\mathtt{Neg}\mathrm (\ulcorner \alpha \urcorner) } 
\wedge \mathtt{Ax}\mathrm (\mathtt{Neg}\mathrm (\ulcorner \alpha 
\urcorner) )] $

	\item  \noindent $[\mathtt{Prf}\mathrm (u) \wedge \mathtt{Fml}\mathrm ((u)_{w}) 
\wedge n = u * 2^{\ulcorner \alpha \urcorner} \wedge 
\mathtt{Gen}\mathrm ((u)_{w},\ulcorner \alpha \urcorner)]
 \text{ and }$
 
 \noindent $  [\mathtt{Prf}\mathrm (u) \wedge \mathtt{Fml}\mathrm 
((u)_{w}) \wedge n = u * 2^{\mathtt{Neg}\mathrm (\ulcorner \alpha 
\urcorner) } \wedge \mathtt{Gen}\mathrm ((u)_{w},\mathtt{Neg}\mathrm 
(\ulcorner \alpha \urcorner) )] $

	\item  $[\mathtt{Prf}\mathrm (u) \wedge \mathtt{Fml}\mathrm ((u)_{z}) 
\wedge \mathtt{Fml}\mathrm ((u)_{w})  \wedge n = u * 2^{\ulcorner 
\alpha \urcorner} \wedge \mathtt{MP}\mathrm 
((u)_{z},(u)_{w},\ulcorner \alpha \urcorner)]
 \text{ and }$
 
 \noindent $
[\mathtt{Prf}\mathrm (u) \wedge \mathtt{Fml}\mathrm ((u)_{z}) \wedge 
\mathtt{Fml}\mathrm ((u)_{w})  \wedge n = u * 2^{\mathtt{Neg}\mathrm 
(\ulcorner \alpha \urcorner)} \wedge \mathtt{MP}\mathrm 
((u)_{z},(u)_{w},\mathtt{Neg}\mathrm (\ulcorner \alpha \urcorner))] $

	\item  $[\mathtt{Prf}\mathrm (u) \wedge  n = u * 2^{\ulcorner \alpha 
\urcorner} \wedge \mathtt{Ax}\mathrm (\ulcorner \alpha \urcorner)]
\text{ and }$

\noindent $  [\mathtt{Prf}\mathrm (u) \wedge  n = u * 
2^{\mathtt{Neg}\mathrm (\ulcorner \alpha \urcorner) } \wedge 
\mathtt{Ax}\mathrm (\mathtt{Neg}\mathrm (\ulcorner \alpha \urcorner) 
)]
$ 
\end{enumerate}
which are all immediately impossible  by the definitions of 
$\mathtt{Ax}\mathrm (y)$,
$\mathtt{Gen}\mathrm (x,y)$ and $\mathtt{MP}\mathrm (x,y,z)$  and 
thence 
by the  definitions of
 the axioms of {\it PA},    the Generalization 
Rule and  Modus Ponens, because no axiom  belongs to {\it PA} 
together with its 
negation and  the two inference rules    preserve
logical validity. 
\end{proof}
We now recall  the definition of characteristic function. If $R$ is a
relation of $n$ arguments, then the characteristic function $C_R$ is 
defined as follows
\begin{equation}\notag
 C_R(x_1,\dots,x_n) =
\begin{cases}
0& \text{if $R(x_1,\dots,x_n)$ is true,}\\
1& \text{if $R(x_1,\dots,x_n)$ is false.}
\end{cases} 
\end{equation} 

Let us call the characteristic functions of $\mathtt{Pf}\mathrm 
(x,v)$, 
$\mathtt{Prf}\mathrm (x)$, $\mathtt{Rf}\mathrm (x,v)$ and  
$\mathtt{Ref}\mathrm (x)$  
respectively $C_{\mathtt{Pf}}$, $C_{\mathtt{Prf}}$, 
$C_{\mathtt{Rf}}$, 
 and $C_{\mathtt{Ref}}$.

A relation $R(x_1,\dots,x_n)$ is said to be primitive recursive 
(recursive) if and only if its
characteristic function $C_R(x_1,\dots,x_n)$ is primitive recursive 
(recursive) (\cite{mendel} 137).
As $\mathtt{Pf}\mathrm (x,v)$, $\mathtt{Prf}\mathrm (x)$,   
$\mathtt{Rf}\mathrm (x,v)$ and $\mathtt{Ref}\mathrm (x)$ are  
primitive recursive 
then also $C_{\mathtt{Pf}}$, $C_{\mathtt{Prf}}$,
 $C_{\mathtt{Rf}}$ and $C_{\mathtt{Ref}}$  are  primitive recursive.

Every recursive function is representable in {\it PA} (\cite{mendel}  
143), thence
 $C_{\mathtt{Pf}}$,  $C_{\mathtt{Prf}}$, 
 $C_{\mathtt{Rf}}$ and $C_{\mathtt{Ref}}$, are
representable in {\it PA}.
We shall  assume   $C_{Pf}$, $C_{Prf}$, 
 $C_{Rf}$ and $C_{Ref}$   to represent respectively 
 $C_{\mathtt{Pf}}$, $C_{\mathtt{Prf}}$, 
 $C_{\mathtt{Rf}}$ and $C_{\mathtt{Ref}}$  in {\it PA}.
\begin{lemma}\label{complete1}
For any formula $\alpha$, and $n$ as the   G\"{o}del number of a proof
in {\it PA} of $\alpha$
$$ 
\vdash _{PA}\: C_{Pf}(\overline{n},\overline{\ulcorner 
\alpha\urcorner})=\overline{0}
 \;\;\wedge\;\; 
C_{Rf}(\overline{n},\overline{\ulcorner \alpha\urcorner})=\overline{1}
$$
\end {lemma} 
\begin{proof}
One can easily see  that the two conjuncts are true: as $n$ is the   
G\"{o}del number of a proof
in {\it PA} of $\alpha$ 
$C_{Pf}(\overline{n},\overline{\ulcorner 
\alpha\urcorner})=\overline{0}$ is true. 
By Lemma (\ref{notboth}) $\mathtt{Rf}\mathrm(n,\ulcorner 
\alpha\urcorner)$ does not hold,
therefore it is true that
 $n$ is not the   G\"{o}del number of a refutation
in {\it PA} of $\alpha$.  
\end{proof}
\begin{lemma}\label{complete2}
For any formula $\alpha$, and $n$ as the   G\"{o}del number of a 
refutation
in {\it PA} of $\alpha$
$$ 
\vdash _{PA}\: C_{Rf}(\overline{n},\overline{\ulcorner 
\alpha\urcorner})=\overline{0}
 \;\;\wedge\;\; 
C_{Pf}(\overline{n},\overline{\ulcorner \alpha\urcorner})=\overline{1}
$$
\end {lemma} 
\begin{proof}
One can easily see  that the two conjuncts are true: as $n$ is the   
G\"{o}del 
number of a refutation
in {\it PA} of $\alpha$ 
$C_{Rf}(\overline{n},\overline{\ulcorner 
\alpha\urcorner})=\overline{0}$ is true. 
By Lemma (\ref{notboth}) $\mathtt{Pf}\mathrm(n,\ulcorner 
\alpha\urcorner)$ does not hold,
therefore it is true that
 $n$ is not the   G\"{o}del number of a proof
in {\it PA} of $\alpha$.  
\end{proof}
\begin{lemma}\label{antidiag1}
For  any formula $\alpha$ 
 
\begin{itemize}
	\item[(i)]\it{ not both }  $$\vdash _{PA}\:  Pf(\overline{n},\overline{\ulcorner 
\alpha\urcorner}) \:\: \vdash _{PA} \: Rf(\overline{n},\overline{\ulcorner 
\alpha\urcorner}),$$
    \item[(ii)] \it{for } $n$  \it{ as the   G\"{o}del number of a refutation
in  PA of } $\alpha$ $$ \vdash _{PA} \: Rf(\overline{n},\overline{\ulcorner \alpha\urcorner}) \Longleftrightarrow 
\neg Pf(\overline{n},\overline{\ulcorner \alpha\urcorner}),$$
	\item[(iii)]\it{for } $n$  \it{ as the   G\"{o}del number of a proof 
	in  PA of } $\alpha  $
$$\vdash _{PA} \:  Pf(\overline{n},\overline{\ulcorner \alpha\urcorner}) 
\Longleftrightarrow \neg Rf(\overline{n},\overline{\ulcorner 
\alpha\urcorner}).$$
\end{itemize}
\end{lemma}
\begin{proof}
(i) Immediately by Lemma (\ref{notboth}) and the definition of being 
expressible which holds for 
both $\mathtt{Pf}\mathrm (x,v)$ and  $\mathtt{Rf}\mathrm (x,v)$ 
(\cite{mendel} 130).

(ii) Let us assume $\vdash _{PA} \: 
Rf(\overline{n},\overline{\ulcorner \alpha\urcorner})$, then
Lemma (\ref{complete2}) yields $\vdash _{PA} \: 
C_{Pf}(\overline{n},\overline{\ulcorner 
\alpha\urcorner})=\overline{1}$.
Hence by definition  $Pf(\overline{n},\overline{\ulcorner 
\alpha\urcorner})$ is false,
consequently $\vdash _{PA} \: \neg 
Pf(\overline{n},\overline{\ulcorner \alpha\urcorner})$.
Conversely let us assume 
$\vdash _{PA} \: \neg Pf(\overline{n},\overline{\ulcorner 
\alpha\urcorner})$
then $Pf(\overline{n},\overline{\ulcorner \alpha\urcorner})$ is false 
and by Lemma
(\ref{complete2}) we attain $\vdash _{PA} \: 
Rf(\overline{n},\overline{\ulcorner \alpha\urcorner})$.

(iii) Let us assume $\vdash _{PA} \: 
Pf(\overline{n},\overline{\ulcorner \alpha\urcorner})$, then
Lemma (\ref{complete1}) yields $\vdash _{PA} \: 
C_{Rf}(\overline{n},\overline{\ulcorner 
\alpha\urcorner})=\overline{1}$.
Hence by definition  $Rf(\overline{n},\overline{\ulcorner 
\alpha\urcorner})$ is false,
consequently $\vdash _{PA} \: \neg 
Rf(\overline{n},\overline{\ulcorner \alpha\urcorner})$.
Conversely let us assume 
$\vdash _{PA} \: \neg Rf(\overline{n},\overline{\ulcorner 
\alpha\urcorner})$
then $Rf(\overline{n},\overline{\ulcorner \alpha\urcorner})$ is false 
and by Lemma
(\ref{complete1}) we attain 
$\vdash _{PA} \: Pf(\overline{n},\overline{\ulcorner 
\alpha\urcorner})$.
\end{proof}

All preceding lemmas were carried out constructively, needlessly  
to assume consistency.
We are now able to consider the consequences yielded by such 
lemmas to the G\"{o}del's  argumentation. 
\begin{itemize}
\item[(a$'$)] Assume $\vdash _{PA} \delta$. 
 Let $r$ be the G\"{o}del number of a proof
in {\it PA} of $\delta$. Then $\mathtt{Pf}$$(r,q)$. Hence, $\vdash 
_{PA} 
Pf(\overline{r},\overline{q})$, that is $\vdash _{PA} 
Pf(\overline{r},\overline{\ulcorner \delta\urcorner})$.
Hence by  Lemma (\ref{antidiag1}) (i) 
$\vdash _{PA} Rf(\overline{r},\overline{\ulcorner \delta\urcorner})$ 
is not admitted, which means that $r$ cannot 
be  the G\"{o}del number of a refutation of $\delta$ 
(indeed Lemma (\ref{complete1}) yields 
$\vdash _{PA} C_{Rf}(\overline{r},\overline{\ulcorner 
\delta\urcorner})=\overline{1}$).
Even though we can have 
$\vdash _{PA} \neg Rf(\overline{r},\overline{\ulcorner 
\delta\urcorner})$,
by (iii) of Lemma (\ref{antidiag1}), then we shall have not 
$\vdash _{PA} \neg Pf(\overline{r},\overline{\ulcorner 
\delta\urcorner})$
 by Lemma (\ref{complete1}) ($\vdash _{PA} 
C_{Pf}(\overline{r},\overline{\ulcorner 
\delta\urcorner})=\overline{0}$).
\item[(b$'$)] Assume  $\vdash _{PA} \neg\delta$. Let $r$
be the G\"{o}del number of a proof in {\it PA} of $\neg \delta$. Then 
$\mathtt{Rf}$$(r,q)$.
Hence $\vdash _{PA} Rf(\overline{r},\overline{q})$ that is
$ \vdash _{PA} Rf(\overline{r},\overline{\ulcorner \delta\urcorner})$.
Hence by  Lemma (\ref{antidiag1})  (i)
$\vdash _{PA} Pf(\overline{r},\overline{\ulcorner \delta\urcorner})$ 
is not admitted.
This means that $r$ cannot be the G\"{o}del number of a proof of 
$\delta$
(in fact, $r$ is  the G\"{o}del number of a refutation of $\delta$,
Lemma (\ref{complete2})  yields
$\vdash _{PA} C_{Pf}(\overline{r},\overline{\ulcorner 
\delta\urcorner})=\overline{1}$).  
Even though we can have
$\vdash _{PA}\neg Pf(\overline{r},\overline{\ulcorner 
\delta\urcorner})$,
by (ii) of Lemma (\ref{antidiag1})  as well, then we shall have not
$\vdash _{PA}\neg Rf(\overline{r},\overline{\ulcorner 
\delta\urcorner})$  
by Lemma (\ref{complete2})  ($\vdash _{PA} 
C_{Rf}(\overline{r},\overline{\ulcorner 
\delta\urcorner})=\overline{0}$).
\end{itemize}
We have thus shown that previous Lemmas  prevent
any accomplishment of (a) and (b) within G\"{o}del's argumentation
\footnote{As regard to  (b), we notice that in (b$'$), by Lemma 
(\ref{notboth}), $\vdash _{PA} Pf(\overline{r},\overline{\ulcorner \delta\urcorner})$ 
is not admitted for every natural number $r$ such that 
$\mathtt{Rf}$$(r,q)$ (i.e. whenever $\vdash _{PA} \neg\delta$).}.
We have then established the following theorem.
\begin{theorem}
By the arithmetization of the refutability predicate
G\"{o}del's  incompleteness does not hold as a theorem  of {\it PA}.
\end{theorem}

\bigskip

\section{Consistency,  Completeness and Decidability}\label{C}
A recursive predicate defines a decidable set, by 
reason that its characteristic function is considered to be 
effectively computable (\cite{mendel} 165, 249).

Let us call $T_{PA}$ the set
of G\"{o}del numbers of theorems of {\it PA} and $R_{PA}$ the set
of G\"{o}del numbers of refutations of {\it PA}.

By the recursiveness of $\mathtt{Pf}\mathrm (x,v)$, 
$C_{\mathtt{Pf}}(x,v)=0$ if $v\in T_{PA}$ 
 and 
$C_{\mathtt{Pf}}(x,v)=1$ if $v\notin T_{PA}$.
By the recursiveness of $\mathtt{Rf}\mathrm (x,v)$, 
$C_{\mathtt{Rf}}(x,v)=0$ if $v\in R_{PA}$ 
 and 
$C_{\mathtt{Rf}}(x,v)=1$ if $v\notin R_{PA}$.

We can than state the following theorem.
\begin{theorem}\label{dec}
 $T_{PA}$  and $R_{PA}$ are decidable sets.
\end{theorem}
It is furthermore 
well-known that
if we have a computable function $f(x_{1} ,\dots,x_{n})$ such that
\begin{equation}\notag
f(x_{1} ,\dots,x_{n}) =
\begin{cases}
0& \text{if $<x_{1} ,\dots,x_{n}>\in S$}\\
1& \text{if $<x_{1} ,\dots,x_{n}>\notin S$}
\end{cases} 
\end{equation} 
(where $S$ is a set of natural number which turns out to be decidable 
just by this definition), 
then the function $g(x_{1} ,\dots,x_{n})$ defined by
\begin{equation}\notag
g(x_{1} ,\dots,x_{n}) =
\begin{cases}
1& \text{if $f(x_{1} ,\dots,x_{n}) =0$}\\
0& \text{if $f(x_{1} ,\dots,x_{n}) =1$}
\end{cases} 
\end{equation}
is effectively computable too. Accordingly the complement of $S$ is 
decidable.
One can easily see that for $f(x_{1} ,\dots,x_{n})$ primitive 
recursive, $g(x_{1} ,\dots,x_{n})$ is primitive recursive too.
Consequently we have 
\begin{equation}\notag
 C_{\mathtt{\neg Prf}}(x) =
\begin{cases}
1& \text{if }C_{\mathtt{Prf}}(x) =0\\
0& \text{if }C_{\mathtt{Prf}}(x) =1
\end{cases}
\end{equation}
\begin{equation}\notag
C_{\mathtt{\neg Ref}}(x) =
\begin{cases}
1& \text{if }C_{\mathtt{Ref}}(x) =0\\
0& \text{if }C_{\mathtt{Ref}}(x) =1
\end{cases} 
\end{equation}
\begin{equation}\notag
 C_{\mathtt{\neg Pf}}(x,v) =
\begin{cases}
1& \text{if }C_{\mathtt{Pf}}(x,v) =0\\
0& \text{if }C_{\mathtt{Pf}}(x,v) =1
\end{cases}
\end{equation}
\begin{equation}\notag
C_{\mathtt{\neg Rf}}(x,v) =
\begin{cases}
1& \text{if }C_{\mathtt{Rf}}(x,v) =0\\
0& \text{if }C_{\mathtt{Rf}}(x,v) =1
\end{cases} 
\end{equation}
where  $\mathtt{\neg Prf}$, $\mathtt{\neg Pf}$, $\mathtt{\neg Ref}$ 
and $\mathtt{\neg Rf}$ are respectively  complementary  of
$\mathtt{ Prf}$, $\mathtt{ Pf}$, $\mathtt{ Ref}$ and $\mathtt{Rf}$.

Let us summarize, $\mathtt{ Prf}$, $\mathtt{ Pf}$, $\mathtt{ Ref}$ and 
$\mathtt{Rf}$ are primitive recursive, then 
$C_{\mathtt{ Prf}}$, $C_{\mathtt{ Pf}}$, $C_{\mathtt{ Ref}}$ and $C_{\mathtt{ Rf}}$
are primitive recursive too. But 
$C_{\mathtt{\neg Prf}}(x)=1- C_{\mathtt{ Prf}}(x)$,
$C_{\mathtt{\neg Pf}}(x,v)=1- C_{\mathtt{ Prf}}(x,v)$,
$C_{\mathtt{\neg Ref}}(x)=1- C_{\mathtt{ Ref}}(x)$,
and
$C_{\mathtt{\neg Rf}}(x,v)=1- C_{\mathtt{ Rf}}(x,v)$,
thence $\mathtt{\neg Prf}$, $\mathtt{\neg Pf}$, 
$\mathtt{\neg Ref}$ and $\mathtt{\neg Rf}$ are primitive recursive too.

We have then  the following statements.
\begin{lemma}\label{Refaut} For every $x$ 
$$\mathtt{Prf}\mathrm (x) 
\text{ if and only if }
\neg\mathtt{Ref}\mathrm (x).$$
\end{lemma}
\begin{proof}
Let us assume 
$\mathtt{Prf}\mathrm (x)$. 
$C_{\mathtt{Prf}}(x) =0$. Hence
 $C_{\mathtt{Prf}}(\mathtt{Neg}(x)) =1$, by 
the 
effective computability of $C_{\mathtt{Prf}}$.
$\mathtt{Prf}(\mathtt{Neg}(x)) $ is false, 
then
$\mathtt{Ref}(x) $ is false. Accordingly, 
$C_{\mathtt{Ref}}(x) =1$. Thus 
$C_{\mathtt{\neg Ref}}(x) =0$ and 
$\neg\mathtt{Ref}\mathrm (x).$

Conversely, let us assume $\neg\mathtt{Ref}\mathrm (x)$.
Then $C_{\mathtt{\neg Ref}}(x) =0$  and
$C_{\mathtt{Ref}}(x) =1$. 
If $\mathtt{Ref}\mathrm (x)$ is false by its 
definition
$\mathtt{Prf}\mathrm (\mathtt{Neg}(x))$ is 
false.
Thus $C_{\mathtt{Prf}}(\mathtt{Neg}(x)) =1$ 
and
$C_{\mathtt{\neg Prf}}(\mathtt{Neg}(x)) =0$.
Consequently
$C_{\mathtt{\neg Prf}}(x) =1$, and
$C_{\mathtt{ Prf}}(x) =0$. Hence
$\mathtt{ Prf}(x)$.
\end{proof}
If we convent to formalize  ``a proof in {\it PA} of $\theta$" with 
 $\theta _{1}\dots\theta _{r}\vdash _{PA}\theta$  then we  have  
 $\vdash _{PA}(\theta _{1}\Rightarrow (\theta _{2}\Rightarrow \dots 
 (\theta _{r}\Rightarrow \theta )\dots))$ (Herbrand, 1930).
 Indeed Lemma (\ref{Refaut}) could be read as follows:
 \emph{ for  $\theta _{1}, \dots , \theta _{r},\theta$  formulas 
in {\it PA} $\mathtt{Prf}\mathrm (\ulcorner(\theta _{1}\Rightarrow 
(\theta _{2}\Rightarrow \dots 
 (\theta _{r}\Rightarrow \theta )\dots))\urcorner)$ if and only if 
$\neg\mathtt{Ref}\mathrm (\ulcorner(\theta _{1}\Rightarrow (\theta 
_{2}\Rightarrow \dots 
 (\theta _{r}\Rightarrow \theta )\dots))\urcorner)$}.

Furthermore, by the recursiveness of $\mathtt{Prf}\mathrm (x)$, 
$C_{\mathtt{Prf}}(x)=0$ if $x\in T_{PA}$ 
 and 
$C_{\mathtt{Prf}}(x)=1$ if $x\notin T_{PA}$.
By the recursiveness of $\mathtt{Ref}\mathrm (x)$
$C_{\mathtt{Ref}}(x)=0$ if $x\in R_{PA}$ 
 and 
$C_{\mathtt{Ref}}(x)=1$ if $x\notin R_{PA}$.
\begin{lemma}\label{autPfautRf} For every $<x,v>$
$$\mathtt{Pf}\mathrm (x,v) 
\text{ if and only if }
\neg\mathtt{Rf}\mathrm (x,v).$$
\end{lemma}
\begin{proof}
Let us assume 
$\mathtt{Pf}\mathrm (x,v)$. 
$C_{\mathtt{Pf}}(x,v) =0$.
Hence 
$C_{\mathtt{Pf}}(x,\mathtt{Neg}(v)) =1$.
Accordingly
$C_{\mathtt{\neg Pf}}(x,\mathtt{Neg}(v)) =0$.
Thence
$C_{\mathtt{\neg Rf}}(x,v) =0$ and
$\mathtt{\neg Rf}(x,v)$.
Conversely, let us assume
$\mathtt{\neg Rf}(x,v)$. 
We have then
$\mathtt{\neg Pf}(x,\mathtt{Neg}(v))$ 
and
$C_{\mathtt{\neg Pf}}(x,\mathtt{Neg}(v)) =0$.
Therefore
$C_{\mathtt{\neg Pf}}(x,v) =1$,  by the 
effective computability of $C_{\mathtt{\neg Pf}}$.
Accordingly
$C_{\mathtt{Pf}}(x,v) =0$ and
$\mathtt{Pf}(x,v)$.
\end{proof}
Indeed Lemma (\ref{autPfautRf}) could  be read as follows:
 \emph{ for  $\theta _{1}, \dots , \theta _{r}, \alpha$ formulas in {\it PA}
$\mathtt{Pf}\mathrm (\ulcorner(\theta _{1}\Rightarrow (\theta 
_{2}\Rightarrow \dots 
 (\theta _{r}\Rightarrow \alpha )\dots)\urcorner,\ulcorner \alpha 
\urcorner)$ if and only if 
$\neg\mathtt{Rf}\mathrm (\ulcorner(\theta _{1}\Rightarrow (\theta 
_{2}\Rightarrow \dots 
 (\theta _{r}\Rightarrow \alpha )\dots)\urcorner,\ulcorner \alpha 
 \urcorner)$}.
\begin{lemma}\label{conscomp}For $m = \ulcorner\alpha\urcorner$ and 
$n = \ulcorner(\theta _{1}\Rightarrow (\theta _{2}\Rightarrow \dots 
 (\theta _{r}\Rightarrow \theta )\dots))\urcorner$
\begin{itemize}
\item[   (i)] $\qquad \qquad m \in T_{PA}\:\text{ iff } \: m \notin R_{PA}, $ 
\item[   (ii)] $\qquad \qquad n \in T_{PA}\:\text{ iff } \: n \notin R_{PA}. $
\end{itemize}
\end{lemma}
\begin{proof}
Immediately (i) by  Lemma (\ref{autPfautRf}), 
 (ii) by Lemma(\ref{Refaut}).
\end{proof}
\begin{theorem}
{\it PA} is consistent; that is, there is no formula $\alpha$ such 
that both
$\alpha$ and $\neg \alpha$ are theorems of {\it PA}.
\end{theorem}
\begin{proof}
Let us assume $m$ to be the G\"{o}del number of a proof of a 
formula $\alpha$ of {\it PA} and $n$ to be the G\"{o}del 
number of a proof of  $\neg \alpha$. Then $n,m \in T_{PA}$.
But since $n$ is the G\"{o}del 
number of a proof of  $\neg \alpha$ we have also $n\in R_{PA}$, 
accordingly
$n$ belongs to both $T_{PA}$ and $R_{PA}$, which is impossible by 
Lemma (\ref{conscomp}).
\end{proof}
\begin{theorem}
{\it PA} is complete; that is for any well formed formula $\alpha$ of 
{\it PA}
either 
$\vdash _{PA} \alpha$ or $\vdash _{PA} \neg \alpha$.
\end{theorem}
\begin{proof}
Let $\alpha$ be a well formed formula of {\it PA}, we can then yield 
by G\"{o}del numbering $m = \ulcorner\alpha\urcorner$. 
 By Lemma (\ref{conscomp}) either 
$m \in T_{PA}$ or $ m \in R_{PA} $. Therefore either 
$\vdash _{PA} \alpha$ or $\vdash _{PA} \neg \alpha$. 
\end{proof}

Let us call {\it PF} the full first-order predicate calculus 
(\cite{mendel} 172).
Let  $T_{PF}$ be then the set of G\"{o}del number of theorems of {\it 
PF}.

\begin{theorem}
$T_{PF}$ is decidable.
\end{theorem}
\begin{proof}
By G\"{o}del Completeness Theorem, a  formula $\alpha$ of {\it PA} is 
provable
in  {\it PA} if and only if $\alpha$ is logically valid, and 
$\alpha$  is provable
in  {\it PF} if and only if $\alpha$ is logically valid. Hence 
$\vdash _{PA} \alpha$ if
and only if $\vdash _{PF} \alpha$. Accordingly, for $n$ as the 
G\"{o}del number of a 
proof of $\alpha$ in {\it PA},
$$n\in T_{PA} \: \text{ iff } \: n\in T_{PF}.$$
Hence, by theorem (\ref{dec}), $T_{PF}$ is decidable.
\end{proof}
Calling our attention to the diagonalization lemma we  note that it 
holds for any formula $\phi$ with only the variable $v$ free. In 
other terms
$\phi$ can be replaced by any formula with only one free variable. 
Let us suppose now that a sentence $\delta$ is a theorem of {\it PA}, i.e. 
$\vdash _{PA} \delta$. For $n$ as the G\"{o}del number of a proof in 
{\it PA} of $\delta$
we have 
$\vdash _{PA} \:  Pf(\overline{n},\overline{\ulcorner 
\delta\urcorner})$.
But for $\phi(v)$ as $\forall x\, Rf(x,v)$  diagonalization lemma 
could have
already yielded
$
\vdash _{\mathit{PA}}\delta\, \iff \,\forall x\, 
Rf(x,\overline{\ulcorner \delta \urcorner})
$
then by biconditional elimination we have
$
\vdash _{\mathit{PA}}\forall x\, Rf(x,\overline{\ulcorner \delta 
\urcorner}).
$
Hence, by Particularization Rule,
$
\vdash _{\mathit{PA}} Rf(\overline{n},\overline{\ulcorner \delta 
\urcorner}),
$
and therefore
$$
\vdash _{\mathit{PA}} Pf(\overline{n},\overline{\ulcorner \delta 
\urcorner}) 
\wedge Rf(\overline{n},\overline{\ulcorner \delta \urcorner}),
$$
which is false 
by reason of   G\"{o}del numbering itself, as proved 
by Lemma (\ref{notboth}) which holds for each natural number $n$. 
By  the tautology $(A \wedge B) \Rightarrow (A \iff B)$  we should 
then have
$
\vdash _{\mathit{PA}} Pf(\overline{n},\overline{\ulcorner \delta 
\urcorner}) 
\iff Rf(\overline{n},\overline{\ulcorner \delta \urcorner}),
$
which openly conflicts with (iii) in Lemma (\ref{antidiag1}).
Finally,
by lemma (\ref{autPfautRf}) it is always the case that for whatever 
formula $\delta$ and $n=\ulcorner(\theta _{1}\Rightarrow (\theta 
_{2}\Rightarrow \dots (\theta _{r}\Rightarrow \delta )\dots)\urcorner$
$$
\vdash _{\mathit{PA}} Pf(\overline{n},\overline{\ulcorner \delta 
\urcorner}) 
\iff \neg Rf(\overline{n},\overline{\ulcorner \delta \urcorner}).
$$
We have thus established that the applicability of the 
diagonalization lemma
to any formula $\phi$ with only the variable $v$ free leads  to
the assertion of a contradiction as a theorem of {\it PA} and for 
this reason {\it PA}  turns out to
be inconsistent. Consequently diagonalization can no longer be 
considered to hold  as an 
equivalence nor replacement theorem. 
We shall have accordingly the following theorem.\begin{theorem}
Diagonalization  does not hold as a lemma in {\it PA}.
\end{theorem} \hfill{$\square$}


\begin{thebibliography}{[ACG95a]}
%
\bibitem{godel1}
G\"{o}del, Kurt. 
\newblock On formally undecidable proposition of Principia 
mathematica and related systems I. 1931. 
\newblock In \emph{Collected Works}, Vol. I Publications 
1929-1936. Oxford University Press, New York, 1986, pp. 145-195.
%
\bibitem{godel2}
G\"{o}del, Kurt. 
\newblock On undecidable proposition of formal
mathematical systems. 1934. 
\newblock In \emph{Collected Works}, Vol. I Publications 1929-1936.
Oxford University Press, New York, 1986, pp. 346-371.
%
\bibitem{hilbert1}
Hilbert, David, Wilhelm Ackermann. 
\newblock \emph{Grundz\"{u}gen der 
theoretischen Logik}. 
\newblock Springer, Berlin, 1928.
%
\bibitem{hilbert2}
Hilbert, David, Paul Bernays.
\newblock \emph{Die Grundlegung der Mathematik I}. 
Springer, Berlin, 1934;
\newblock \emph{Die Grundlegung der Mathematik II}. 
Springer, Berlin, 1939.
%
\bibitem{mendel}
Mendelson, Elliott.
\newblock  \emph{Introduction to Mathematical Logic} Third 
Edition. 
\newblock The Wads-worth \& Brooks, Pacific Grove, California, 1987.
\end{thebibliography}
\end{document}